\documentclass{article}[11pt,titlepage]

\newtheorem{theorem}{Theorem}

\newcommand\nc{\newcommand}
\nc{\od}[2]{\frac{d#1}{d#2}}
\nc{\be}{\begin{equation}}
\nc{\ee}{\end{equation}}
\nc{\bd}{\begin{displaymath}}
\nc{\ed}{\end{displaymath}}
\nc{\BEA}{\begin{eqnarray}}
\nc{\EEA}{\end{eqnarray}}
\nc{\p}{\partial}

\usepackage{graphics}
\usepackage{amsmath}
\usepackage{amsfonts}
\usepackage{amssymb}

\begin{document}
\title{A one-dimensional model for the interaction between cell-to-cell adhesion and chemotactic signalling}
\author{K. Anguige\thanks{e-mail: kmpa@hotmail.com} ,\\ {\em Wolfgang Pauli Institute,}\\ {\em Fakult\"at f\"ur Mathematik, Universit\"at Wien}\\ {\em Nordbergstra\ss{}e 15, 1090 Wien, Austria.}}

\date{Revised January 12th, 2011}

\maketitle

\begin{abstract}
We develop and analyse a discrete, one-dimensional model of cell motility which incorporates the effects of volume filling, cell-to-cell adhesion and chemotaxis. The formal continuum limit of the model is a nonlinear generalisation of the parabolic-elliptic Keller-Segel equations, with a diffusivity which can become negative if the adhesion coefficient is large. The consequent ill-posedness results in the appearance of spatial oscillations and the development of plateaus in numerical solutions of the underlying discrete model. A global-existence result is obtained for the continuum equations in the case of favourable parameter values and data, and a steady-state analysis which, amongst other things, accounts for high-adhesion plateaus is carried out. For ill-posed cases, a singular Stefan-problem formulation of the continuum limit is written down and solved numerically, and the numerical solutions are compared with those of the original discrete model. 
\end{abstract}

\section{Introduction}

Of late, there has been considerable interest in formulating continuum models for cell structures generated by cell-to-cell adhesion, the motivation being to facilitate mathematical analysis and efficient numerical simulation of processes such as {\em de novo} blood-vessel synthesis ({\em i.e.} vasculogenesis) and cancer invasion, for example. Two recent attempts in this direction were made, respectively, by Armstrong {\em et al.} \cite{xi}, who began with a nonlocal integro-differential equation in which the kernel is integrated over a given cell-sensing radius, and by Anguige \& Schmeiser \cite{AS}, who wrote down a simple random-walk model accounting for adhesion, diffusion and volume filling. For both approaches, it turns out that the limiting macroscopic model (a nonlinear parabolic equation) can be ill posed if the adhesion is sufficiently strong, which leads to interesting pattern-forming behaviour in solutions of the underlying microscopic models, but which also makes mathematical analysis rather more difficult than one would like \cite{A,AS}.

Our intention in this paper is to extend the modelling and analysis of \cite{A,AS} in a rather obvious way, namely, by factoring in the directed response of cells to an extracellular chemical gradient ({\em i.e.} {\em chemotaxis}), and then examining the resulting interaction between such {\em long-range} signalling and the {\em short-range} signalling of cell-to-cell adhesion.

Following the approach previously adopted in \cite{AS}, our 1-d model for cell adhesion, diffusion and chemotaxis will take the form of the random-walk system
\be
\frac{\partial \rho_i}{\partial t}=\mathcal{T}^+_{i-1}\rho_{i-1}+\mathcal{T}^-_{i+1}\rho_{i+1}-(\mathcal{T}^+_i+\mathcal{T}^-_i)\rho_i~,\label{walk}
\ee
for the approximate cell densities $\rho_i\in[0,1]$, on a uniform grid of points $x_i=ih$, the quantities $\mathcal{T}^{\pm}_i$ being the transitional probabilities per unit time of a one-step jump from $i$ to $i\pm 1$.

Taking inspiration from \cite{AS}, as well as \cite{v}, the scaled transitional probabilities are chosen to be
\be
\mathcal{T}^{\pm}_i = \frac{1}{h^2}(1-\rho_{i\pm 1})(1-\alpha\rho_{i\mp 1})\left(1+\frac{\chi_0}{2}(S_{i\pm 1}-S_i)\right)~,\label{T_i} 
\ee
where $\alpha\in[0,1]$ is the adhesion coefficient, $\chi_0\in[0,\infty)$ the chemotactic sensitivity, and $S_i$ the concentration of chemoattractant at the point $x_i$. Here, the first factor in parentheses models volume filling, the second models adhesion, and the third chemotaxis. The justification for the adhesion factor and the range of $\alpha$ is that, for example, the presence of a particle to the right should reduce the probability of a particle jumping to the left: $\alpha<0$ would correspond to repulsion, and $\alpha>1$ would allow the ${T}^{\pm}_i$ to go negative.

The expression (\ref{T_i}) includes as special cases (and is a simple combination of) both our previous model for adhesion/diffusion ($\chi_0=0$) \cite{AS}, and the model for linear diffusion, chemotaxis and volume filling ($\alpha=0$) presented in \cite{v}. 

Upon taking the continuum limit of (\ref{walk})-(\ref{T_i}) by writing, for example, $\rho_{i\pm 1}=\rho(x_i\pm h)$ and $S_{i\pm 1}=S(x_i\pm h)$, and taking Taylor expansions in powers of $h$, one ends up with the advection-diffusion equation 
\be 
\frac{\p\rho}{\p t} = \frac{\p}{\p x}\left(D(\rho)\frac{\p\rho}{\p x}-\chi(\rho)\rho\frac{\p S}{\p x}\right),\label{rho}
\ee
where the diffusivity is given by
\be
D(\rho) = 3\alpha\left(\rho-\frac{2}{3}\right)^2 + 1 -\frac{4}{3}\alpha,\label{D}
\ee
and the chemotactic-sensitivity function by
\be
\chi(\rho) = \chi_0(1-\rho)(1-\alpha\rho)\label{chi}.
\ee

The continuum equation for the chemoattractant is taken to be the usual non-dimensionalised quasi-steady-state equation
\be
\Delta S = S-\rho.\label{S}
\ee

Equations (\ref{rho})-(\ref{S}) constitute a nonlinear generalisation of the classical Keller-Segel chemotaxis model, depending on just two parameters, $\chi_0$ and $\alpha$. They are to be solved on the domain $\Omega=(0,L)$, subject to the homogeneous Neumann conditions $\rho_x=S_x=0$ at $x=0,L$.

Note that the presence of the volume-filling term in $\chi(\rho)$ entails that (\ref{rho}) has the following maximum principle: 
\be
0\leq\rho(0)\leq 1\Rightarrow 0\leq\rho(t)\leq 1. 
\ee

Moreover, (\ref{S}) implies a maximum principle for $S$:
\be
0\leq\rho(t)\leq 1\Rightarrow 0\leq S(t)\leq 1.
\ee

As in \cite{AS}, however, (\ref{rho}) can be ill-posed if $\alpha>\frac{3}{4}$, since in that case there is an interval of values of $\rho$, centred around $\rho=\frac{2}{3}$, for which the diffusivity is negative. Explicitly, we have from (\ref{D}) that $D(\rho)<0$ whenever
\be
\rho\in I_{\alpha}:= \left(\frac{2\alpha-\sqrt{\alpha(4\alpha-3)}}{3\alpha},\frac{2\alpha+\sqrt{\alpha(4\alpha-3)}}{3\alpha}\right) = (\rho^{\flat},\rho^{\sharp}),\label{rhointerval}
\ee
and $D(\rho)\geq 0$ otherwise. Note that the width of $I_{\alpha}$ increases as $\alpha\nearrow 1$, and that $I_1=\left(\frac{1}{3},1\right)$.

Equation (\ref{rho}), therefore, is certainly ill posed if the initial density profile hits $I_{\alpha}$, and we saw in \cite{AS} that this ill-posedness is related to the presence of oscillations and the formation of plateaus in solutions of (\ref{walk}) in the special case $\chi_0=0$. As a consequence of chemotactic aggregation, however, we might also expect (\ref{rho})-(\ref{S}) to be ill posed for small initial densities, provided $\chi_0$ is large enough and $\alpha>\frac{3}{4}$.

What interests us most in this paper, then, is possible singular pattern-forming behaviour in the case $\alpha>\frac{3}{4}$ and $\chi_0>0$. This turns out to be rather difficult to investigate analytically, but considerable insight may nevertheless be gained with the aid of numerical simulations.

The paper is organised as follows.

In Section 2 we carry out a steady-state analysis of (\ref{rho})-(\ref{S}), and this is followed in Section 3 by some global-existence results for favourable parameter values and initial data.

In Section 4 we report on numerical simulations of (\ref{rho})-(\ref{S}), which show that singular ({\em i.e.} sharp-edged) aggregation patterns can be generated by small data, provided the chemotactic sensitivity and the adhesion coefficient are chosen large enough.

Finally, and in analogy with \cite{A}, we consider the idea that a Stefan-problem-type framework, in which the density is allowed to jump across the unstable region $I_{\alpha}$, might be an appropriate way of treating (\ref{rho})-(\ref{S}) as the limit of (\ref{walk}) in ill-posed cases. Although such problems seem to be analytically intractable at present, there is a sense in which solving them numerically may nevertheless be more efficient than discretising the Neumann problem for (\ref{rho})-(\ref{S}) directly, since this (the direct approach) requires a very fine mesh to properly resolve the singular behaviour observed near $I_{\alpha}$; simulations obtained using both of these methods are compared and contrasted in Section 5.

\section{Steady-state analysis}

We now show that essentially the same techniques employed in the case of linear diffusion (see, {\em e.g.}, \cite{v}), along with a comparison-principle argument, can be used to investigate steady states of (\ref{rho})-(\ref{S}). For $\alpha>\frac{3}{4}$, it is also possible, as in the case $\chi_0=0$ \cite{AS}, to construct discontinuous weak solutions in which $\rho$ has finitely many jumps across the unstable region $I_{\alpha}$.

\subsection{Linear stability of uniform steady states}

Linearising (\ref{rho}), (\ref{D}), (\ref{chi}), (\ref{S}) around a uniform steady state $\rho=S=\bar{\rho}$, and inserting the ansatz $\rho=e^{\lambda t}e^{ik\pi x/L},~S=Ae^{\lambda t}e^{ik\pi x/L}$, gives $A=L^2/(L^2+k^2\pi^2)$ and the dispersion relation
\be
\lambda = \frac{k^2\pi^2}{L^2}\left(-D(\bar{\rho}) + \frac{L^2\chi(\bar{\rho})\bar{\rho}}{L^2+k^2\pi^2}\right),
\ee
for the growth rate $\lambda$ and the wave number $k$.

Thus, if $\chi_0$ is so small that $\chi(\bar{\rho})\bar{\rho}<D(\bar{\rho})$, then $\lambda<0~\forall k$, and the uniform steady state is  unconditionally stable. If, on the other hand, $\chi(\bar{\rho})\bar{\rho}>D(\bar{\rho})$, which necessarily occurs when $\alpha>\frac{3}{4}$ and $\bar{\rho}\in I_{\alpha}$, for example, then $\lambda>0$ for small wave numbers, and the uniform steady state is thus unstable to long-wavelength perturbations. 

It is also elementary to show that the dominant wavemode is determined by
\be
\left(\frac{k\pi}{L}\right)^2 = \sqrt{\frac{\chi(\bar{\rho})\bar{\rho}}{D(\bar{\rho})}} - 1.
\ee

\subsection{Global $L^{\infty}$-stability of uniform steady states for $\alpha<\frac{3}{4}$ and $\chi_0$ small}

For a given mass $M:=\|\rho\|_{L^1}$, (\ref{rho})-(\ref{S}) always has the uniform solution $(\rho,S)=(\bar{\rho},\bar{\rho})$, where $\bar{\rho}=M/L$. One expects such a solution to be a global attractor provided $\alpha<\frac{3}{4}$ and $\chi_0$ is sufficiently small, and this is the content of the following theorem.

\begin{theorem}
If 
\be
\chi_0\min\left\{1,\sqrt{L/2}\right\} \max_{0\leq\rho\leq 1}(1-\rho)(1-\alpha\rho)\rho<(1-\frac{4}{3}\alpha),\label{stab_cond}
\ee
then any smooth, global solution pair $(\rho,S)$ of (\ref{rho})-(\ref{S}) satisfies $\|\rho-\bar{\rho}\|_{\infty}(t)\leq c_1\|\rho-\bar{\rho}\|_2(0)e^{-c_2t}$ and $\|S-\bar{\rho}\|_{H^2}(t)\leq \|\rho-\bar{\rho}\|_2(0)e^{-c_2t}$, for some positive constants $c_1,~c_2$.
\end{theorem}

{\em Proof.}~~~~First of all, in  (\ref{rho}), we subtract $\bar{\rho}$ from $\rho$, multiply through by $\rho-\bar{\rho}$, and integrate by parts to get
\be
\frac{1}{2}\frac{d}{dt}\|\rho-\bar{\rho}\|_2^2 = -\int_0^L D(\rho)((\rho-\bar{\rho})_x)^2~dx + \chi_0\int_0^L (\rho-\bar{\rho})_x(1-\rho)(1-\alpha\rho)\rho S_x~dx\label{L^2_est}.
\ee

Next, note that squaring (\ref{S}) gives, with the aid of an integration by parts,
\be
\int_0^L (S-\bar{\rho})^2 + 2((S-\bar{\rho})_x)^2 + ((S-\bar{\rho})_{xx})^2~dx = \int_0^L(\rho-\bar{\rho})^2~dx,\label{H2}
\ee
and that differentiating (\ref{S}) and carrying out the same procedure gives
\be
\int_0^L ((S-\bar{\rho})_x)^2 + 2((S-\bar{\rho})_{xx})^2 + ((S-\bar{\rho})_{xxx})^2~dx = \int_0^L((\rho-\bar{\rho})_x)^2~dx,\label{H3}
\ee
since the Neumann conditions kill the boundary terms in both cases.

In particular, (\ref{H2}) and (\ref{H3}), together with the Poincar\'{e} inequality, imply that
\be
\|(S-\bar{\rho})_x\|_2\leq\min\left\{1,\sqrt{L/2}\right\}\|(\rho-\bar{\rho})_x\|_2.
\ee
Hence, condition (\ref{stab_cond}) and a further application of the Poincar\'{e} inequality imply that
\BEA
\frac{1}{2}\frac{d}{dt}\|\rho-\bar{\rho}\|_2^2 & \leq & -\epsilon\|(\rho-\bar{\rho})_x\|_2^2  \nonumber \\
& \leq & -\frac{\epsilon}{L}\|\rho-\bar{\rho}\|_2^2,\label{Poincare}
\EEA
for some $\epsilon>0$, and therefore also
\be
\|\rho-\bar{\rho}\|_2(t)\leq\|\rho-\bar{\rho}\|_2(0)e^{-\epsilon t/L},\label{L^2_decay}
\ee
and, by (\ref{H2}),
\be
\|S-\bar{\rho}\|_{H^2}\leq\|\rho-\bar{\rho}\|_2(0)e^{-\epsilon t/L}.\label{S_H2_est}
\ee

Next, for a given $n\in\mathbb{N}$, we integrate the first of (\ref{Poincare}) from $t=n$ to $n+1$, thus obtaining
\be
\int_n^{n+1}\limits\|(\rho-\bar{\rho})_x\|_2^2~dt\leq C\|\rho-\bar{\rho}\|_2^2(0)e^{-\epsilon n/L},
\ee
which implies that $\forall n\in\mathbb{N}, \exists~\zeta^n\in[n,n+1]$ such that $\|(\rho-\bar{\rho})_x\|_2(\zeta^n)\leq C\|\rho-\bar{\rho}\|_2(0)e^{-\epsilon n/L}$. Therefore, by Sobolev imbedding and the $L^2$-decay estimate (\ref{L^2_decay}) (or Poincar\'{e} again),
\be
\|\rho-\bar{\rho}\|_{\infty}(\zeta^n)\leq C\|\rho-\bar{\rho}\|^2_2(0)e^{-2\epsilon n/L}.
\ee

From here, we proceed with a comparison argument, which will demonstrate that $\rho-\bar{\rho}$ can grow (pointwise) at worst exponentially between the times $\zeta^n$, which of course satisfy $|\zeta^n-\zeta^{n+1}|\leq 2$. 

Expanding (\ref{rho}), and using (\ref{S}) to substitute for the Laplacian, we see that $u:=\rho-\bar{\rho}$ satisfies 
\be
\frac{\p u}{\p t} = \frac{\p}{\p x}\left(D(u+\bar{\rho})\frac{\p u}{\p x}\right) - \frac{\p}{\p\rho}\chi(\rho)\rho\frac{\p u}{\p x}\frac{\p S}{\p x} - \chi(\rho)\rho(-u + (S-\bar{\rho})).
\ee

Thus, by (\ref{S_H2_est}), we have that $u$ is a subsolution of the following problem, for $t\in[\zeta^n,\zeta^{n+1})$:
\be
\frac{\p w}{\p t} = \frac{\p}{\p x}\left(D(x,t)\frac{\p w}{\p x}\right) + A(x,t)\frac{\p w}{\p x} + B(x,t)w + C\|\rho-\bar{\rho}\|_2(0)e^{-\epsilon t/L},\label{comp}
\ee
subject to 
\be
w_x(0,t)=w_x(L,t)=0,\qquad w(x,\zeta^n) = C\|\rho-\bar{\rho}\|_2(0)e^{-\epsilon\zeta^n/L},\label{comp_data}
\ee
where $B(x,t)$ is a bounded function, say $|B(x,t)|\leq C_0~\forall~x,t$, and $A(x,t)$,~$D(x,t)$ are smooth functions, with $D(x,t)>0$.

It is easy to check that $\hat{w}:=Ce^{-\epsilon\zeta^n/L}\|\rho-\bar{\rho}\|_2(0)e^{C_1(t-\zeta^n)}$ is a supersolution of (\ref{comp})-(\ref{comp_data}) provided $C_1>C_0$ is chosen sufficiently large.

Subsolutions can be constructed similarly, and we therefore arrive at
\be
\|u\|_{\infty}(t)\leq C\|\rho-\bar{\rho}\|_2(0)e^{-\epsilon t/L},\label{u_inf}
\ee
$\forall t>0$, as required. $\square$
\\~\\
The argument used in the proof of Theorem 1 can now be bootstrapped to obtain
\begin{theorem}
The uniform steady state, $(\rho,S)=(\bar{\rho},\bar{S})$, of (\ref{rho})-(\ref{S}) is a local $L^{\infty}$-attractor provided
\be
\min\left(1,\sqrt{L/2}\right)\chi(\bar{\rho})\bar{\rho}<D(\bar{\rho}).\label{boot}
\ee
\end{theorem}
 
{\em Proof.}~~~~Let $E(\rho):=D(\rho)-\min\left(1,\sqrt{L/2}\right)\chi(\rho)\rho$, and set $\delta:=D(\bar{\rho})-\min\left(1,\sqrt{L/2}\right)\chi(\bar{\rho})\bar{\rho}>0$. By the proof of Theorem 1, it is enough to show that $E(\rho)$ cannot hit $\frac{\delta}{2}$ in finite time along a solution trajectory, provided $\|\rho_0-\bar{\rho}\|_{\infty}$ is chosen small enough. Assume the contrary, and let $t^\ast$ be the first time at which $E(\rho)=\frac{\delta}{2}$ for a given initial datum $\rho_0$. Then, for $t<t^\ast$, one has an inequality of the form (\ref{u_inf}), and hence, if $\|\rho_0-\bar{\rho}\|_{\infty}$ was chosen small enough, $E(\rho)>\frac{3\delta}{4}$ for $t<t^\ast$, which is a contradiction. $\square$

Note that condition (\ref{boot}) is consistent with the linear-stability analysis of Section 2.1, and that this condition holds for $\chi_0$ sufficiently small, provided either $\alpha<\frac{3}{4}$ or both $\alpha>\frac{3}{4}$ and $\bar{\rho}\notin I_{\alpha}$.

\subsection{Non-trivial steady states; dynamical systems}

\subsubsection{The case $\alpha<\frac{3}{4}$}
We are looking for smooth solutions of
\BEA
D(\rho)\frac{\p\rho}{\p x} - \chi(\rho)\rho\frac{\p S}{\p x} &  =  & 0, \label{steady_rho}\\
S_{xx} - S + \rho  & = & 0,\label{steady_S}
\EEA
subject to the boundary conditions $S_x=\rho_x=0$ at $x=0,L$.

First of all, equation (\ref{steady_rho}) can be integrated to get
\be
G(\rho) = S - C,\label{G}
\ee
where $G$ is a primitive for $D(\rho)/(\chi(\rho)\rho)$ ($G$ strictly increasing since $\alpha<\frac{3}{4}$), and $C$ is a constant of integration.

Using (\ref{steady_S}), this gives us
\be
G(\rho)_{xx} = G(\rho) - \rho + C,\label{G_DS}
\ee
or, in terms of $\sigma:=G(\rho)$,
\be
\sigma_{xx} = \sigma - G^{-1}(\sigma) + C,\label{Ham}
\ee
which is a Hamiltonian dynamical system, whose critical points are necessarily saddles or centres.

By considering the shape of $G(\rho)-\rho$, it is easy to see that (\ref{Ham}) has either one or three critical points, depending on $\alpha, \chi_0$ and the choice of $C$. If $\chi_0$ is small, then there is just one critical point, and the only possible steady states satisfying the boundary conditions are uniform solutions. If, on the other hand, $\chi_0$ is sufficiently large, then the function $G(\rho)-\rho$ has two extrema, and for an interval of values of $C$ there are three critical points of (\ref{Ham}), comprising two saddles and a nonlinear centre inbetween. Solutions of the Neumann problem for (\ref{Ham}) are then realised as half-orbits around the centre, or integer multiples thereof.

If we have obtained a solution, $\rho$, of (\ref{G_DS}), and if $S$ is then determined by the Neumann problem for (\ref{steady_S}), then it is easy to see that the original steady-state equation (\ref{steady_rho}) also holds. Indeed, using (\ref{steady_S}) to substitute for the linear $\rho$-term in the rhs of (\ref{G_DS}) leads to the energy equality
\be
\|G(\rho)-S+C\|_{H^1}=0,
\ee
which gives the desired result, since $G'(\rho) = D(\rho)/(\chi(\rho)\rho)$.

By inspecting $G(\rho)$, it is also easy to see that there must be sufficient (and also not too much) mass for such solutions to exist.

Moreover, using the fact that the Hamiltonian $F(\sigma):=\int^{\sigma}\hat{\sigma} - G^{-1}(\hat{\sigma})+C~d\hat{\sigma}$ increases towards the centre (denoted by $(\sigma,\sigma_x)=(\sigma_c,0)$), one can, by integrating the Hamiltonian equation, easily show that the length of half-orbits is minimised as $\sigma\rightarrow\sigma_c$, and that the minimum such length, $L^{\ast}$, is given by half the period of linearised simple harmonic motion, $\sigma_{xx}=F''(\sigma_c)\sigma$, about $\sigma_c$. Thus, for a given $C$ there is a minimum value of the domain length which allows for non-trivial steady-state solutions, and it is given by
\be
L^{\ast} = \pi\sqrt{\frac{D(\rho_c)}{(\chi(\rho_c)\rho_c-D(\rho_c))}},
\ee
where $\rho_c=G^{-1}(\sigma_c)$.

In fact, there is a critical curve in $(\alpha,\chi_0)$-space which divides the parameter region in which there is the possibility of three  critical points from that in which there can only be a single one. It is obtained by solving $F(\rho)=F'(\rho)=0$, where $F(\rho):=D(\rho)-\chi(\rho)\rho$, and is plotted in Figure 1 as the boundary between regions (ii) and (iii).

\begin{figure}

\centering

\resizebox{4.5in}{4.5in}{\includegraphics{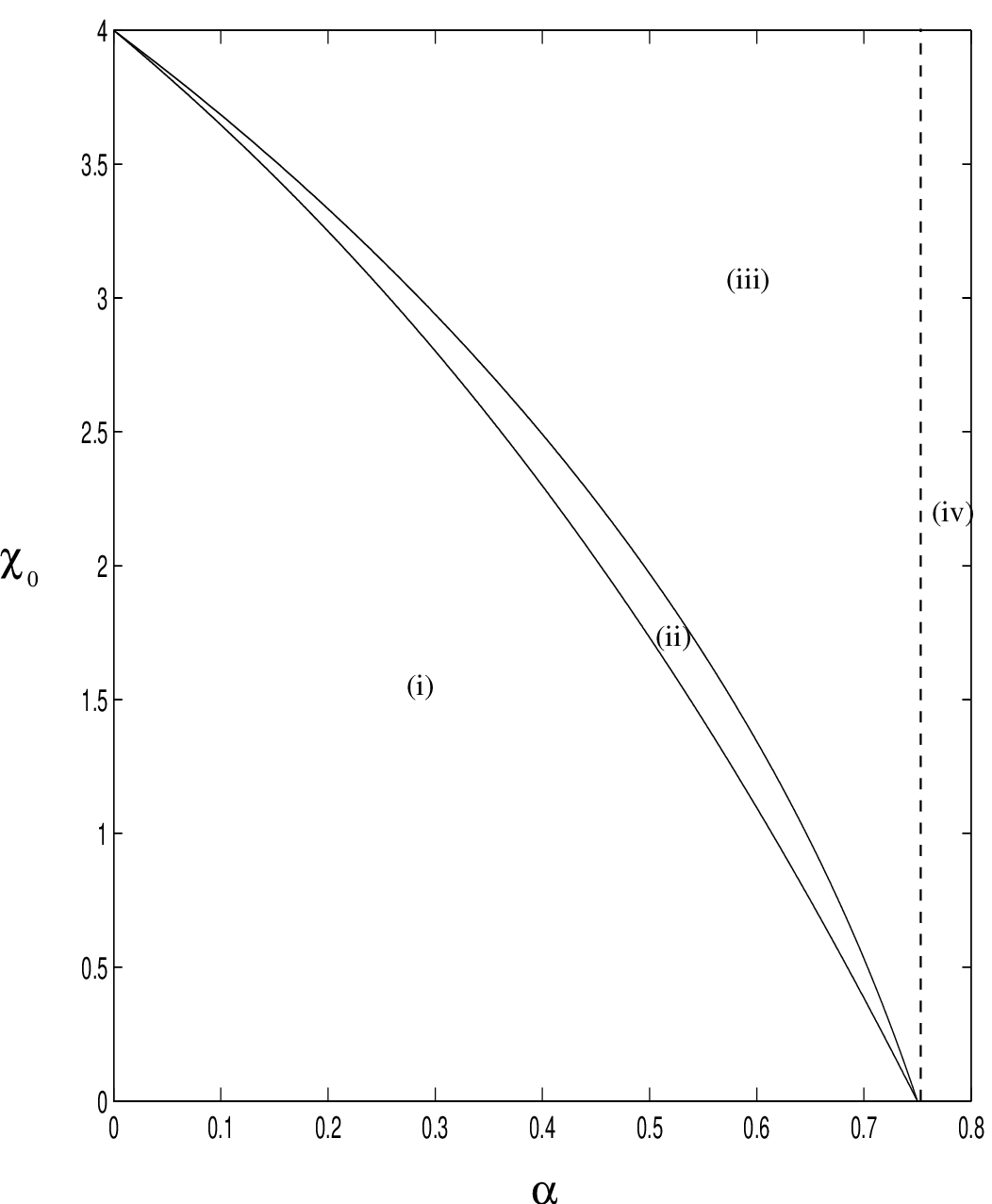}}

\caption{$(\alpha,\chi_0)$-space.}

\end{figure}

A fuller explanation of Figure 1 is as follows. Region (i) is where there are no non-uniform steady states, and where uniform steady states have been proved to be nonlinearly stable, no matter how large $L$ is (see Theorem 1). In region (ii) there are still no non-uniform steady states, and uniform steady states are always linearly stable; nonlinear stability has only been proved for small enough L (Theorem 1 again). In region (iii), non-uniform steady states become possible, for appropriate $C$ and large enough $L$, and uniform states can lose their stability. In region (iv), the diffusivity can turn negative, and one may look for both smooth and discontinuous steady-state solutions, as we now discuss.

\subsubsection{The case $\alpha>\frac{3}{4}$; weak solutions}

For $\alpha>\frac{3}{4}$, non-uniform smooth solutions which miss $I_{\alpha}$ are possible for all choices of $\chi_0>0$, provided $L$ is sufficiently large. For small $\chi_0$ these solutions have to be small-amplitude oscillations just outside the unstable region $I_{\alpha}$, while for large $\chi_0$ they have to be small-amplitude oscillations near $\rho=0$ or $\rho=1$. These conclusions are again reached by considering the form of $G(\rho)$.

We can also look for weak solutions $(\rho,S)$ of (\ref{steady_rho})-(\ref{steady_S}) such that $\rho$ contains finitely many jumps across $I_{\alpha}$, but is smooth and satisfies (\ref{steady_rho}) elsewhere. In that case, $S$, as determined globally by (\ref{steady_S}), will still remain $H^2$ ($\Rightarrow C^1$) smooth, by elliptic regularity.

A natural weak formulation of (\ref{steady_rho}) would be to look for $\rho\in L^2$ satisfying
\be
\int_0^L \phi_{xx}K(\rho) + \phi_x\chi(\rho)\rho S_x~dx = 0,\qquad\forall\phi\in C_0^{\infty}(0,L),\label{weak_steady}
\ee
where $K(\rho)$ is a primitive for $D(\rho)$.

If we suppose that $\rho$ has discontinuities at $x=s_i$, $i=1,...,n$, then the jump conditions for a solution of (\ref{weak_steady}) satisfying the Neumann conditions at $x=0,L$ are calculated to be
\be
(\chi(\rho)\rho S_x - D(\rho)\rho_x)(s_i^{\pm}) = 0\quad\textrm{and}\quad K(\rho(s_i^-))=K(\rho(s_i^+)).\label{jump_cond}
\ee

The first of (\ref{jump_cond}) implies, by the continuity of $S_x$, that the density gradients on either side of a discontinuity are coupled {\em via}
\be
\frac{D(\rho)}{\chi(\rho)\rho}\rho_x(s_i^-) = \frac{D(\rho)}{\chi(\rho)\rho}\rho_x(s_i^+).\label{jump_2}
\ee

Candidate weak solutions can be constructed by solving, on each interval $(s_i,s_{i+1})$, the dynamical system
\be
G(\rho)_{xx} = G(\rho) - \rho + C_i,\label{Gi_DS}
\ee
for suitable $C_i\in\mathbb{R}$ (see below), such that either $\rho<\rho^{\flat}$ or $\rho>\rho^{\sharp}$ in each phase, and patching together the orbits to get a discontinuous density $\rho$ on $[0,L]$. Once this is done, $S$ is globally determined by solving the Neumann problem for (\ref{steady_S}), as already mentioned.

For a pair $(\rho,S)$ determined in this way, we can recover equation (\ref{steady_rho}) on each interval, as in Section 2.3.1. Thus, substituting (\ref{steady_S}) into (\ref{Gi_DS}) for each $i$, we arrive at
\be
\sum_{i=1}^n\left\{\|G(\rho) - S + C_i\|^2_{H^1((s_i,s_{i+1}))} - [(G-S+C_i)(G-S+C_i)_x]_i\right\} = 0,
\ee 
where $[\cdot]_i$ denotes the leap at $x=s_i$.

It follows that $G(\rho)=S-C_i$ in each phase if, for example, the endpoints of adjacent orbits are chosen so that both the flux condition (\ref{jump_2}) and the jump condition $[(G(\rho) + C_i)]_i=0$ hold. A global solution, $\rho$, constructed in this way is also a weak solution in the sense of (\ref{weak_steady}) iff neighbouring $C_i$ are chosen so that $[K(\rho)]_i=0~\forall i$. 
\\~\\
{\em Remark.} Steady states of the Stefan-problem formulation for (\ref{rho})-(\ref{S}), to be introduced in Section 5, are just special cases of such weak solutions for which the endpoints $\rho(s_i^{\pm})$ take on a particular pair of $\alpha$-dependent values; see Section 5.1 for details.

\section{Existence and uniqueness results}

\subsection{Global existence for $\alpha<\frac{3}{4}$}

In the case $\alpha=0$, a global-existence theorem for (\ref{rho}), (\ref{S}) follows directly from the ideas of \cite{HP}. We now extend this result to cover all $\alpha\in\left[0,\frac{3}{4}\right)$.
\begin{theorem}
Given smooth initial data $(\rho_0,S_0)$ satisfying $0\leq\rho_0\leq 1, 0\leq S_0\leq 1$, the system (\ref{rho}), (\ref{D}), (\ref{chi}), (\ref{S}) has a unique smooth global solution, $(\rho,S)$, satisfying $0\leq\rho\leq 1,~0\leq S\leq 1$, provided $\alpha<\frac{3}{4}$.
\end{theorem}

{\em Proof.}~~~~Since smooth data satisfying the Neumann boundary condition at $x=0,L$ can be reflected about $x=0$ to give $H^3$ data on the circle $S^1$, and since (\ref{rho})-(\ref{S}) are invariant under the transformation $x\mapsto -x$, it is enough to prove our existence theorem for $H^3$ data on $S^1$.

For this, first note that (\ref{S}) can be used to write 

\be
S = (1-\Delta)^{-1}\rho,
\ee
and that $(1-\Delta)^{-1}$ is a bounded operator from $L^2(S^1)$ to $H^2(S^1)$. This can now be substituted into (\ref{rho}), thus reducing our problem to the single nonlocal diffusion equation

\be 
\frac{\p\rho}{\p t} = \frac{\p}{\p x}\left(D(\rho)\frac{\p\rho}{\p x}-\chi(\rho)\rho(1-\Delta)^{-1}\frac{\p\rho}{\p x}\right).\label{rho'}
\ee

This equation can be solved using mollifiers, in line with the treatment of quasi-linear parabolic equations presented in Section 15.7 of \cite{x}. Specifically, given $H^s$ data $\rho_0$, we first of all introduce the equation

\be
\frac{\p\rho_{\epsilon}}{\p t} = J_{\epsilon}\frac{\p}{\p x}\left(D(J_{\epsilon}\rho_{\epsilon})\frac{\p J_{\epsilon}\rho_{\epsilon}}{\p x}-\chi(J_{\epsilon}\rho_{\epsilon})J_{\epsilon}\rho_{\epsilon}(1-\Delta)^{-1}\frac{\p J_{\epsilon}\rho_{\epsilon}}{\p x}\right),\label{moll}
\ee
with initial data $J_{\epsilon}\rho_0$, where $J_{\epsilon}$ is a Friedrichs mollifier, which, in particular, is self-adjoint in $L^2$, commutes with $\frac{\p}{\p x}$, is uniformly bounded in each $C^k$ and $H^s$ for $\epsilon<1$, and satisfies $\|J_{\epsilon}u-u\|_{H^s}\rightarrow 0$ as $\epsilon\rightarrow 0$, for any $u\in H^s$. 

Equation (\ref{moll}) can be regarded as an ODE in $H^s$ for $\rho_{\epsilon}$, such that the right-hand side is Lipschitz continuous in this quantity. Local existence and uniqueness of solutions follows by Picard's Theorem, as applied to Banach spaces \cite{ix}.

One now aims to get high-order uniform Sobolev bounds on the family $\rho_{\epsilon}$, and thus to obtain a solution of (\ref{rho'}) as the limit of a convergent subsequence $\rho_{\epsilon_k}$.

The main point to emphasise here is that, by $L^2\rightarrow H^2$ smoothing, a source term such as $(1-\Delta)^{-1}\rho_x$ is, in terms of Sobolev norms, `better than' $\rho$, and that there is consequently no obstacle to proceeding as in Chapter 15 of \cite{x} to derive the required estimates.

To get the $L^2$ estimate, for example, we multiply (\ref{moll}) by $\rho_{\epsilon}$, and integrate by parts, resulting in
\BEA
\frac{1}{2}\frac{d}{dt}\|\rho_\epsilon\|_2^2 & = &\int_{S^1}-D(J_{\epsilon}\rho_{\epsilon})|\p_xJ_{\epsilon}\rho_{\epsilon}|^2 + (\p_xJ_{\epsilon}\rho_{\epsilon})\chi(J_{\epsilon}\rho_{\epsilon})J_{\epsilon}\rho_{\epsilon}(1-\Delta)^{-1}\p_x J_{\epsilon}\rho_{\epsilon}~dx \nonumber \\
& \leq & -\left(1-\frac{4}{3}\alpha\right)\|\p_x J_{\epsilon}\rho_{\epsilon}\|_2^2 + C(\|\rho_{\epsilon}\|_{\infty})\|\p_x J_{\epsilon}\rho_{\epsilon}\|_2\|J_{\epsilon}\rho_{\epsilon}\|_2 \nonumber \\
& \leq & C(\|\rho_{\epsilon}\|_{\infty})\|\rho_{\epsilon}\|_2^2,
\EEA
where we used the elementary inequality
\be
AB\leq C_0A^2+B^2/4C_0 \label{elem}
\ee
in the last line.

Following \cite{x} closely, one can, by repeatedly differentiating (\ref{moll}) and using Moser estimates for product and composite functions, arrive at the family of estimates
\be
\frac{d}{dt}\|\rho_{\epsilon}(t)\|^2_{H^s}\leq C_l(\|\rho_{\epsilon}(t)\|_{C^2})(\|\rho_{\epsilon}(t)\|^2_{H^s}+1).\label{rho_est}
\ee 

By comparing (\ref{rho_est}) with its associated ODE, and using compact Sobolev imbedding into $C^2$ for $s>\frac{5}{2}$, this gives, as in Lemma 7.1, Theorems 7.2, 7.4, Chapter 15 of \cite{x}, a sequence $\rho_{\epsilon_k}\rightarrow\rho\in C([0,T],C^2(S^1))\cap C^{\infty}((0,T)\times S^1)$, such that $\rho$ solves (\ref{rho'}) on some time interval $[0,T]$.

Next, Theorems 8.3, 9.6 and 9.10 of (\cite{x}, Ch.15), along with the Maximum Principle, imply that a continuation criterion for (\ref{rho'}) is that the chemotaxis term be bounded in $L^p$, some $p>\frac{n}{2}$, on finite time intervals, $n$ being the space dimension. Thus, global existence of solutions will be established if we can show that the quantity
\be
\left\|\frac{\p}{\p x}\left((1-\rho)(1-\alpha\rho)\rho\frac{\p S}{\p x}\right)\right\|_{L^1(S^1)}
\ee
is bounded on each $[0,T]$. For this, it is sufficient to show that $\|\frac{\p\rho}{\p x}\|_{L^1}$ is bounded on each $[0,T]$, since $S_x$ and $S_{xx}$ are {\em a priori } uniformly bounded for all $t$, by the Maximum Principle applied to (\ref{S}).

We proceed as in Lemma 2.2 of \cite{viii}, and define an approximation of the sign function by
\be
\sigma_{\delta}(z) = \sigma(z/\delta),\quad\textrm{for}\quad 0<\delta\ll 1,
\ee
with $\sigma$ a smooth and increasing function such that $\sigma(0)=0$ and
\be
\sigma(z)=\textrm{sgn}(z)\quad\textrm{for}\quad |z|>z_0,
\ee
some $z_0>0$.

Then, setting $\textrm{abs}_{\delta}(z):=\int_0^z\sigma_{\delta}(\xi)~d\xi$, we get $\textrm{abs}_{\delta}(z)\rightarrow |z|$ as $\delta\rightarrow 0$, uniformly in $z$. Also note that 
\be
\sigma_{\delta}'(z)=\left\{
\begin{array}{ccl}
0 & , & z>z_0\delta\\
O(\frac{1}{\delta}) & , & z\leq z_0\delta.
\end{array}
\right.\label{sigma'}
\ee

Next, we differentiate (\ref{rho}) w.r.t. $x$, and mutiply by $\sigma_{\delta}(\p_x\rho)$, which, upon integration, leads to 
\be
\frac{d}{dt}\int_0^1\textrm{abs}_{\delta}(\p_x\rho)~dx = Q_{\delta} + \int_0^1\sigma_{\delta}(\p_x\rho)\p^2_{xx}(D(\rho)\rho_x)~dx,\label{abs_delta}
\ee
where $Q_{\delta}$ comes from the chemotaxis term, and can be treated as in \cite{viii} since the argument there does not depend on the specific form of $\chi(\rho)$, resulting in
\be
\lim_{\delta\rightarrow 0}|Q_{\delta}| \leq C_1 + C_2\int_0^1 |\p_x\rho|~dx.
\ee
For the remaining term on the rhs of (\ref{abs_delta}), we get, writing $v=\rho_x$, and integrating by parts,
\be
\begin{array}{cl}
& \int_0^1\sigma_{\delta}'(v)(-D(\rho)v_x^2 - D'(\rho)v_xv^2)~dx\\
 & \\
\leq & \int_0^1\sigma_{\delta}'(v)(-D(\rho)v_x^2 + |D'(\rho)||v_x|v^2)~dx\\
 & \\
\leq & C\int_0^1\sigma_{\delta}'(v)v^4~dx\\
 & \\
\leq & C\delta^3,
\end{array}
\ee
where we have used the condition $\sigma_{\delta}(0)=0$ to kill the boundary terms, along with $\sigma_{\delta}'\geq 0$, equations (\ref{sigma'}), (\ref{elem}), and the fact that $D(\rho)\geq 1-\frac{4}{3}\alpha>0$.

Thus, taking the limit $\delta\rightarrow 0$ in (\ref{abs_delta}), we see that
\be
\frac{d}{dt}\int_0^1|\p_x\rho|~dx\leq C_1 + C_2\int_0^1|\p_x\rho|~dx.
\ee

An application of Gronwall's inequality now shows that the continuation criterion is satisfied.

Finally, for uniqueness, suppose that $(\rho,S)$ and $(\hat{\rho},\hat{S})$ are two solution pairs for (\ref{rho})-(\ref{S}) with the same data. Subtract the equation satisfied by $\hat{S}$ from that satisfied by $S$, square both sides and integrate by parts to obtain
\be
\int_0^L (S-\hat{S})^2 + 2(S_x-\hat{S}_x)^2 + (S_{xx}-\hat{S}_{xx})^2~dx = \int_0^L (\rho-\hat{\rho})^2~dx, \label{Sdiff}
\ee
and subtract the equation satisfied by $\hat{\rho}$ from that satisfied by $\rho$ to get
\be
\frac{\p}{\p t}(\rho-\hat{\rho}) = \frac{\p}{\p x}\left(D(\rho)\frac{\p}{\p x}(\rho-\hat{\rho}) + \frac{\p\hat{\rho}}{\p x}(D(\rho)-D(\hat{\rho})) - \chi(\rho)\rho\frac{\p}{\p x}(S-\hat{S}) - (\chi(\rho)\rho-\chi(\hat{\rho})\hat{\rho})\frac{\p\hat{S}}{\p x}\right).\label{unique}
\ee
Now multiply (\ref{unique}) by $\rho-\hat{\rho}$, and use integration by parts, the mean-value theorem, (\ref{elem}) and (\ref{Sdiff}) to get
\be
\frac{1}{2}\frac{d}{dt}\|\rho-\hat{\rho}\|_2^2\leq C\|\rho-\hat{\rho}\|_2^2.
\ee
Gronwall's inequality implies that $\rho=\hat{\rho}$, and then (\ref{Sdiff}) gives $S=\hat{S}$, as required. $\square$
\\~\\
One might also expect global existence of smooth solutions for $\alpha>\frac{3}{4}$, provided the initial density profile were uniformly outside $I_{\alpha}$ (either above or below) and $\chi_0$ were small enough. This can actually be proved for data close enough to a uniform steady state, as we now demonstrate.

\subsection{Global existence for $\alpha>\frac{3}{4}$, $\chi_0$ small, and small data}

It turns out that we can use the comparison argument of Theorem 1 to obtain a global-existence theorem for (\ref{rho})-(\ref{S}) when $\alpha>\frac{3}{4}$, provided $\chi_0$ is small enough and the initial density profile, $\rho_0(x)$, is sufficiently far from the unstable interval $I_{\alpha}$. We can obtain results for both of the cases $\rho_0(x)<\rho^{\flat}~\forall x\in[0,L]$ and $\rho_0(x)>\rho^{\sharp}~\forall x\in[0,L]$, but for clarity we will simply concentrate on the case $\rho_0(x)<\rho^{\flat}$ in what follows. 

First note that a local-in-time solution is guaranteed by the $J_{\epsilon}$-method used previously, and that for global existence we merely need to prevent the solution from hitting $\rho^{\flat}$.

Thus, let $\bar{\rho}$ be given, and pick $\rho_{\delta_1}>\bar{\rho}$ such that $D(\rho)\geq\delta_1$ for $\rho\leq\rho_{\delta_1}$. Next, pick an initial datum $\rho_0(x)$ such that $\textrm{avg}(\rho_0):= \frac{1}{L}\int_0^L\rho_0(x)dx =\bar{\rho}$, and such that $\max\rho_0(x)\leq\rho_{\delta_1}$, and introduce a smooth, modified diffusivity $D^*(\rho)$ which is equal to $D(\rho)$ for $\rho\leq\rho_{\delta_1}$, and which is greater than $\frac{1}{2}\delta_1$ for $\rho\geq\rho_{\delta_1}$.

Equation (\ref{rho})$^*$ is then defined by replacing $D(\rho)$ with $D^*(\rho)$ in the right-hand side of (\ref{rho}). Global existence of a solution $(\rho,S)$ to (\ref{rho})$^*$, (\ref{chi}), (\ref{S}) follows by previous arguments, and it remains to show that $\rho\leq\rho_{\delta_1}~\forall t$, provided $\chi_0$ and $\|\rho_0-\bar{\rho}\|_{\infty}$ are chosen small enough.

For this, let
\be
\delta_2 = \chi_0\max_{0\leq\rho\leq\rho^{\flat}}(1-\rho)(1-\alpha\rho)\rho,
\ee
and define
\be
\epsilon=\frac{1}{2}\delta_1-\delta_2.
\ee

If $\chi_0$ is chosen so small that $\epsilon>0$, then, in the same way as for $\alpha<\frac{3}{4}$, we obtain the $L^2$-decay estimate
\be
\|\rho-\bar{\rho}\|_2(t)\leq\|\rho_0-\bar{\rho}\|_2 \, e^{-\epsilon t/L},
\ee
along with a sequence $\zeta^n\in[n,n+1]$, such that 
\be
\|\rho-\bar{\rho}\|_{\infty}(\zeta^n)\leq C(\delta_1,\chi_0)\|\rho_0-\bar{\rho}\|_2\label{infty_2}.
\ee
An $H^2$-estimate of the form (\ref{S_H2_est}) again holds, whereby we emphasise that the right-hand side is $O(\|\rho_0-\bar{\rho}\|_2)$, and it therefore follows from (\ref{infty_2}) and the comparison argument of Theorem 1 that
\be
\|\rho-\bar{\rho}\|_{\infty}(t)\leq C(\delta_1,\chi_0)\|\rho_0-\bar{\rho}\|_2.
\ee
Thus, $\rho\leq\rho_{\delta_1}~\forall t$ if $\|\rho_0-\bar{\rho}\|_{\infty}$ is chosen small enough, and consequently $(\rho,S(\rho))$ solves (\ref{rho})-(\ref{S}) for all time.

Given a globally existing solution, convergence to the uniform steady state can be also be proved by the method of Theorem 1, provided $\chi_0$ is sufficiently small.

The case $\rho>\rho^{\sharp}$ is handled analogously, and we therefore have, in summary, 

\begin{theorem}
Given a smooth initial datum $\rho_0(x)$ satisfying either $\rho_0(x)<\rho^{\flat}~\forall x\in[0,L]$ or $\rho_0(x)>\rho^{\sharp}~\forall x\in[0,L]$, and letting $\bar{\rho}:=\textrm{avg}(\rho_0)$, equations (\ref{rho})-(\ref{S}) have a unique, local-in-time solution $(\rho,S)$, which continues to exist globally if both $\|\rho_0-\bar{\rho}\|_{\infty}$ and $\chi_0$ are small. The smallness required of $\chi_0$ depends on $\bar{\rho}$ and $\|\rho_0-\bar{\rho}\|_{\infty}$. Furthermore, given any globally existing smooth solution which misses $I_{\alpha}$, long-time exponential $L^{\infty}$-convergence to the uniform steady state holds, provided $\chi_0$ is sufficiently small.
\end{theorem}

It should be noted that this result does not rule out the possibility of a globally existing solution with $\alpha>\frac{3}{4}$ approaching some non-uniform steady state outside $I_{\alpha}$ as $t\rightarrow\infty$; in this regard, see the simulations in the next section.

\section{Numerics for the discrete model}

In this section we numerically solve the Neumann problem for (\ref{rho})-(\ref{S}) in the high-adhesion regime by means of a finite-difference scheme on a uniform spatial grid of $n$ mesh points, $x_i$, a distance $h$ apart.

The discretisation of the diffusion term in (\ref{rho}) is obtained by setting $\chi_0=0$ in the right-hand side of (\ref{walk}), while the chemotaxis term is discretised by means of the simple, $O(h)$-accurate upwinding scheme
\be
\frac{\p}{\p x}\left(\chi(\rho)\rho\frac{\p S}{\p x}\right)(x_i)\approx\left\{
\begin{array}{ccc}
L_{i+\frac{1}{2}}-L_{i-\frac{1}{2}}  &: & v_i\geq 0 \\
L_{i+\frac{3}{2}}-L_{i+\frac{1}{2}} &: & v_i < 0
\end{array}
\right.,
\ee
where
\be
L_{i+\frac{1}{2}} = \frac{1}{2}\rho_i(1-\rho_i)(1-\alpha\rho_i)(S_{i+1}-S_{i-1})/h^2,\qquad v_i = S_{i+1}-S_{i-1}.
\ee
Thus, in order to obtain numerical stability in well-posed regions, we have chosen a method of lines which is slightly different from (but $O(h)$-consistent with) the original discrete model (\ref{walk}) - if one simply uses the whole of (\ref{walk}), then, upon time integration, the lack of upwinding leads to spurious oscillations, even when $\alpha<\frac{3}{4}$. In contrast, the oscillations caused by our specific discrete adhesion model when $\alpha>\frac{3}{4}$ are to be thought of as fundamental, particular to the model, and what we are really interested in.

To complete the finite-difference scheme, we solve the elliptic equation (\ref{S}) at each time step {\em via} the usual discrete Laplacian, together with Matlab matrix inversion. The solution is updated by means of a semi-implicit time discretisation which is in the spirit of \cite{AS}. That is to say, obvious linear factors of $\rho_i$ in the rhs of the $i$-th $\rho$-equation are evaluated at the new time, rather than the old.

In what follows, we are particularly interested in observing how our numerical solutions change in the vicinity of $I_{\alpha}$ as $n$ increases, since this is essentially the same as asking in what sense (\ref{rho}) is the continuum limit of (\ref{walk}) when $\alpha>\frac{3}{4}$.

\subsection{Singular aggregation patterns}

We begin by choosing $\alpha=0.95, \chi_0=16$, and discretising a small initial density profile on a domain of length $L=8$, using a grid of $n=400$ spatial points. Evolving this data with our numerical scheme results in the sequence of snapshots displayed in Figure 2. 

Clearly, the effect of chemotaxis is to draw the solution towards $I_{\alpha}$, and once it has penetrated sufficiently far into the unstable region, a small number of fine oscillations quickly develop. Subsequently, mass is sucked into the central oscillatory region
{\em via} a combination of chemotaxis and backward diffusion, while positive diffusion flattens out the density profile on either side. Eventually, after a slow process of coarsening in which the fine oscillations disappear, we are left with a single, sharp-edged plateau, which presumably represents a steady-state weak solution of (\ref{rho})-(\ref{S}), as constructed in Section 2.3.2.

One further point to note here is that the values between which $\rho$ jumps at the plateau edge are very close to those observed in \cite{AS} for the case $\chi_0=0$, where we also saw oscillations, as well as plateau formation through coarsening. In that paper the fact that the jump values appeared to depend only on $\alpha$, and not on the initial data, was explained by the existence of a unique heteroclinic cycle for an $O(h^2)$ modified equation derived from (\ref{walk}). We claim that the same argument goes through in the case $\chi_0>0$, since the chemotactic terms produce a higher-order correction due to the smoothing properties of (\ref{S}). The relevant saddle-point values for $\alpha=0.95$ are, in notation that will be used again in Section 5, $(\rho_1,\rho_2)=(0.055,0.99)$ (see Figure 8 of \cite{AS}).

\begin{figure}

\centering

\resizebox{5in}{6in}{\includegraphics{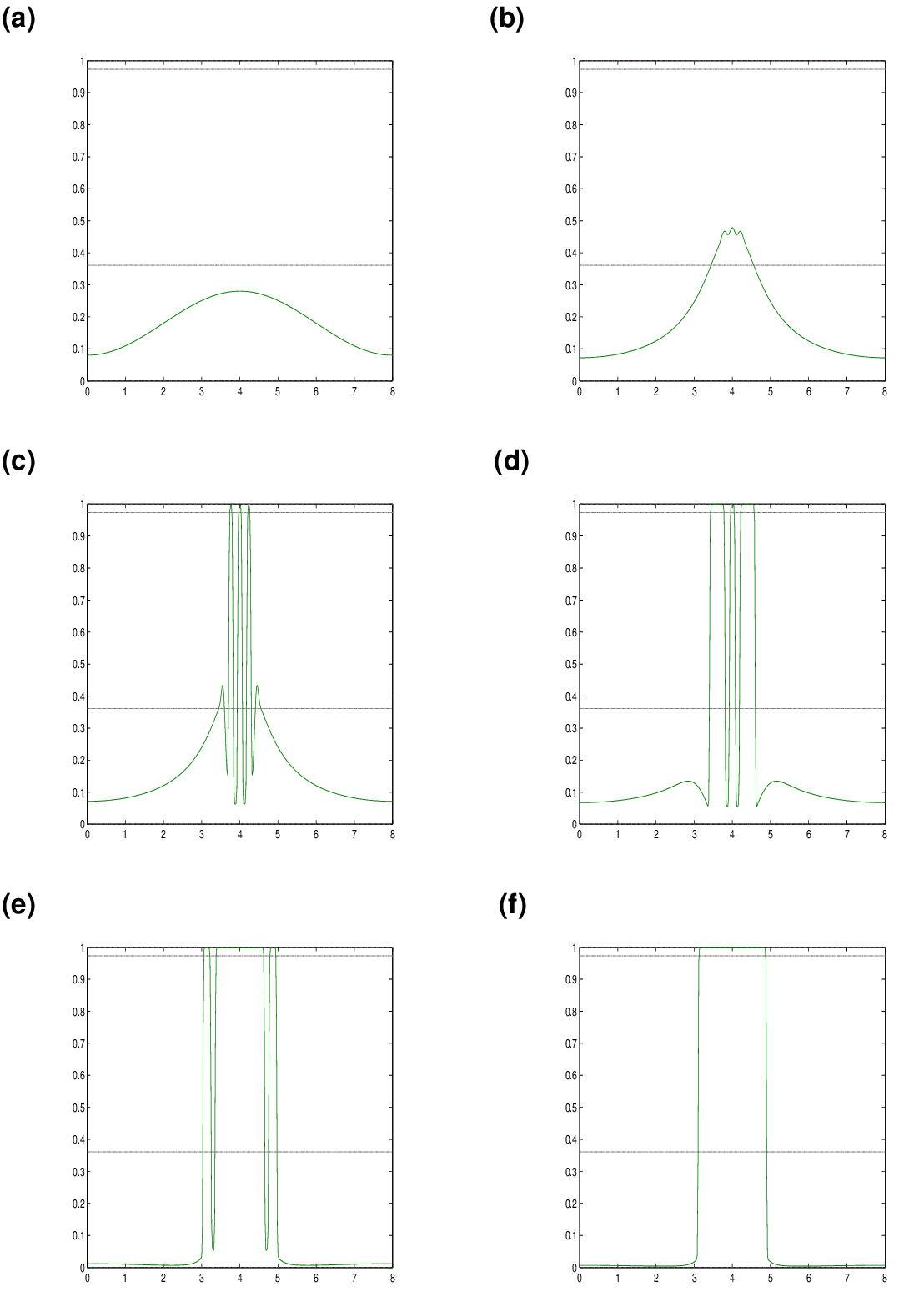}}

\caption{Evolution of a small initial density profile, using $n=400$ spatial points. Data shown at (a) $t=0$, (b) $t=1.3$, (c) $t=1.38$, (d) $t=1.8$, (e) $t=6$, (f) $t=7$. The parameter values are $\alpha=0.95$, $\chi_0=16$, $L=8$, and the boundaries of the unstable region $I_{\alpha}$ are marked with dotted lines.}

\end{figure}

Next, we repeat the simulation of Figure 2, but this time with $n=800$ spatial points; the resulting snapshots are depicted in Figure 3.
Notable differences with respect to Figure 2 are that the oscillations appear a little earlier, and at a slightly lower density level, and that the jumps levels at large times appear even closer to $\rho_1$ and $\rho_2$.

\begin{figure}

\centering

\resizebox{5in}{6in}{\includegraphics{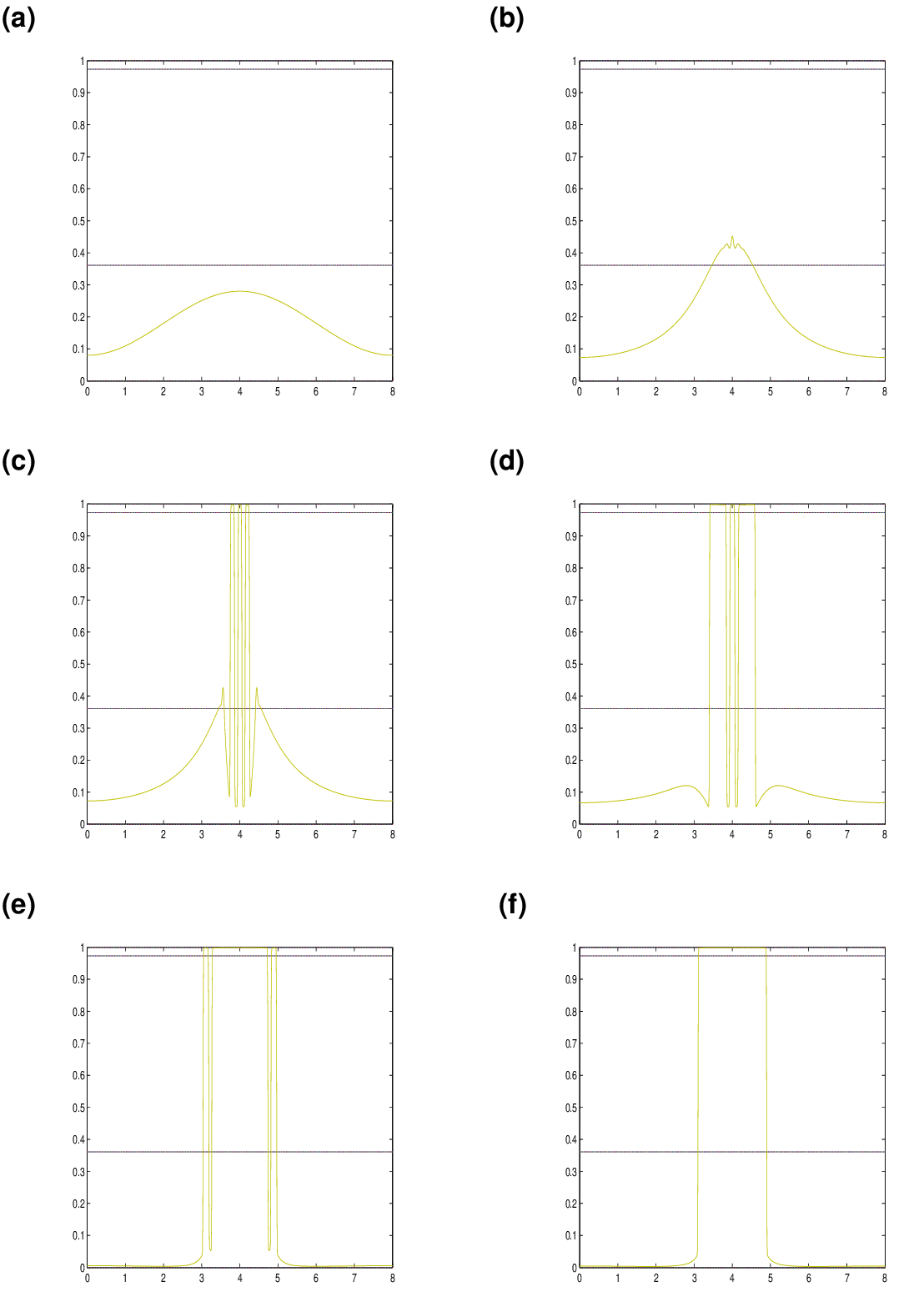}}

\caption{Evolution of the same data as in Figure 2 using $n=800$ spatial points. Data shown at (a) $t=0$, (b) $t=1.1675$, (c) $t=1.255$, (d) $t=1.805$, (e) $t=7.2$, (f) $t=7.6$. The parameter values are $\alpha=0.95$, $\chi_0=16$, $L=8$.}

\end{figure}

Finally, we repeat the simulation once more using $n=1200$ spatial points; the solution is plotted in Figure 4. Note that the level at which oscillations arise is now even earlier, and closer to the lower end of $I_{\alpha}$, but that the rate of convergence to (presumably) $\rho^{\flat}$ is exceedingly slow w.r.t $n$. The jump levels have closed in even further on $\rho_1$ and $\rho_2$. 

\begin{figure}

\centering

\resizebox{5in}{6in}{\includegraphics{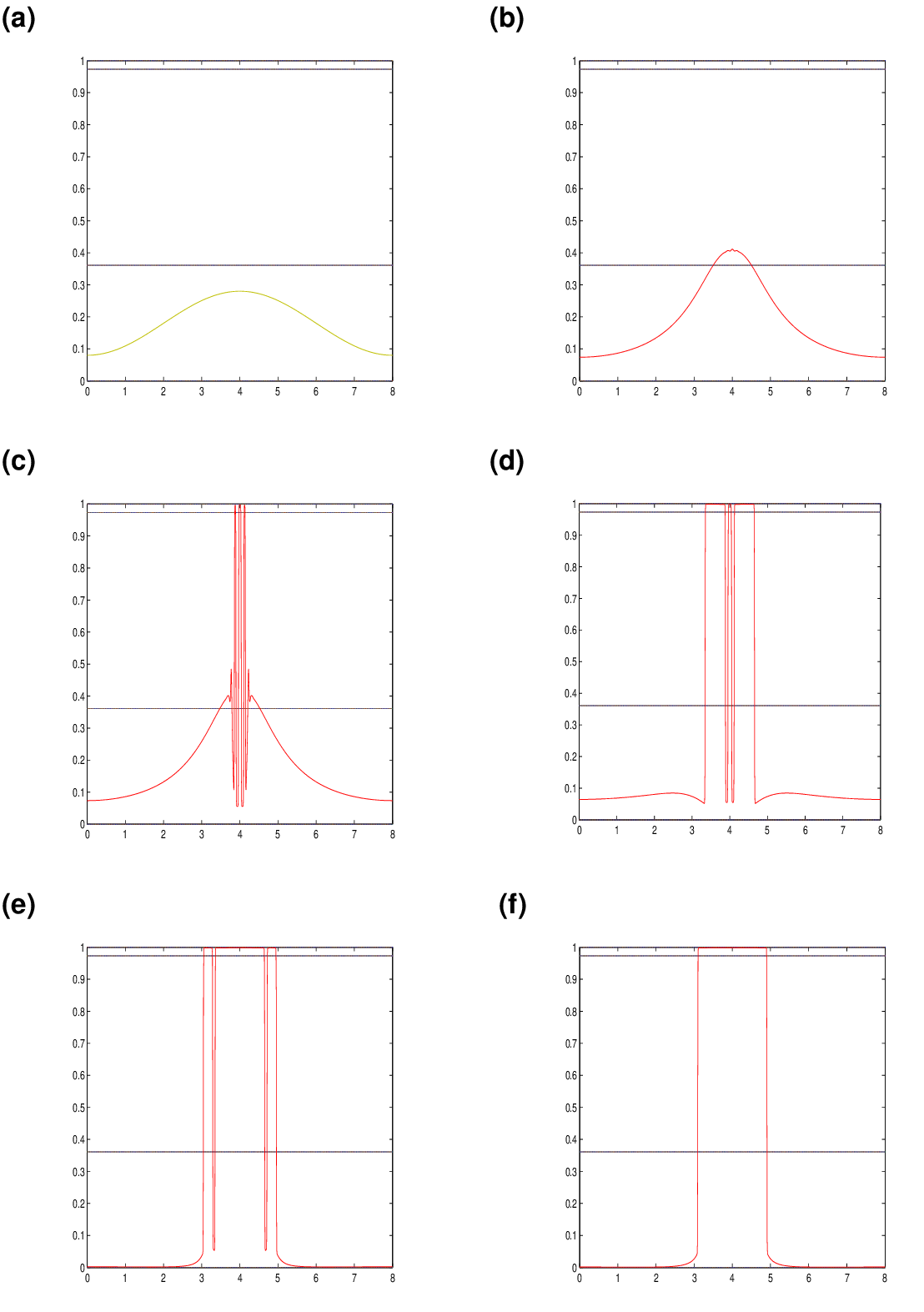}}

\caption{Evolution of the same data as in Figure 2, using $n=1200$ spatial points. Data shown at (a) $t=0$, (b) $t=1.0911$, (c) $t=1.1356$, (d) $t=2$, (e) $t=8.6778$, (f) $t=9.344$. The parameter values are $\alpha=0.95$, $\chi_0=16$, $L=8$.}

\end{figure}

\subsection{Macroscopic coarsening}

Putting aside microscopic oscillations (which are related to ill-posedness), for the moment, a phenomenon exhibited by other, well-posed chemotaxis models is that of macroscopic coarsening, whereby a large aggregation region attracts a smaller one (see, {\em e.g.}, \cite{v,viii}). Such behaviour is in fact also exhibited by our model (\ref{rho})-(\ref{S}), as evidenced by the simulation of Figure 5, in which a wide plateau region absorbs a much narrower neighbour, resulting in a (quasi?) steady state at large times.

This brings us to another phenomenon associated with chemotaxis equations subject to the Neumann condition, which is that a single, asymmetrical plateau will tend to move (perhaps very slowly) towards the boundary as $t\rightarrow\infty$. Unfortunately, our numerical code is not accurate enough to say definitively which way (if any) the plateau in Figure 5d is moving; one might hazard a guess that the plateau remains where it is, due to the fact that there is essentially a Dirichlet condition on either side of the jump locations.

\begin{figure}

\centering

\resizebox{5in}{5in}{\includegraphics{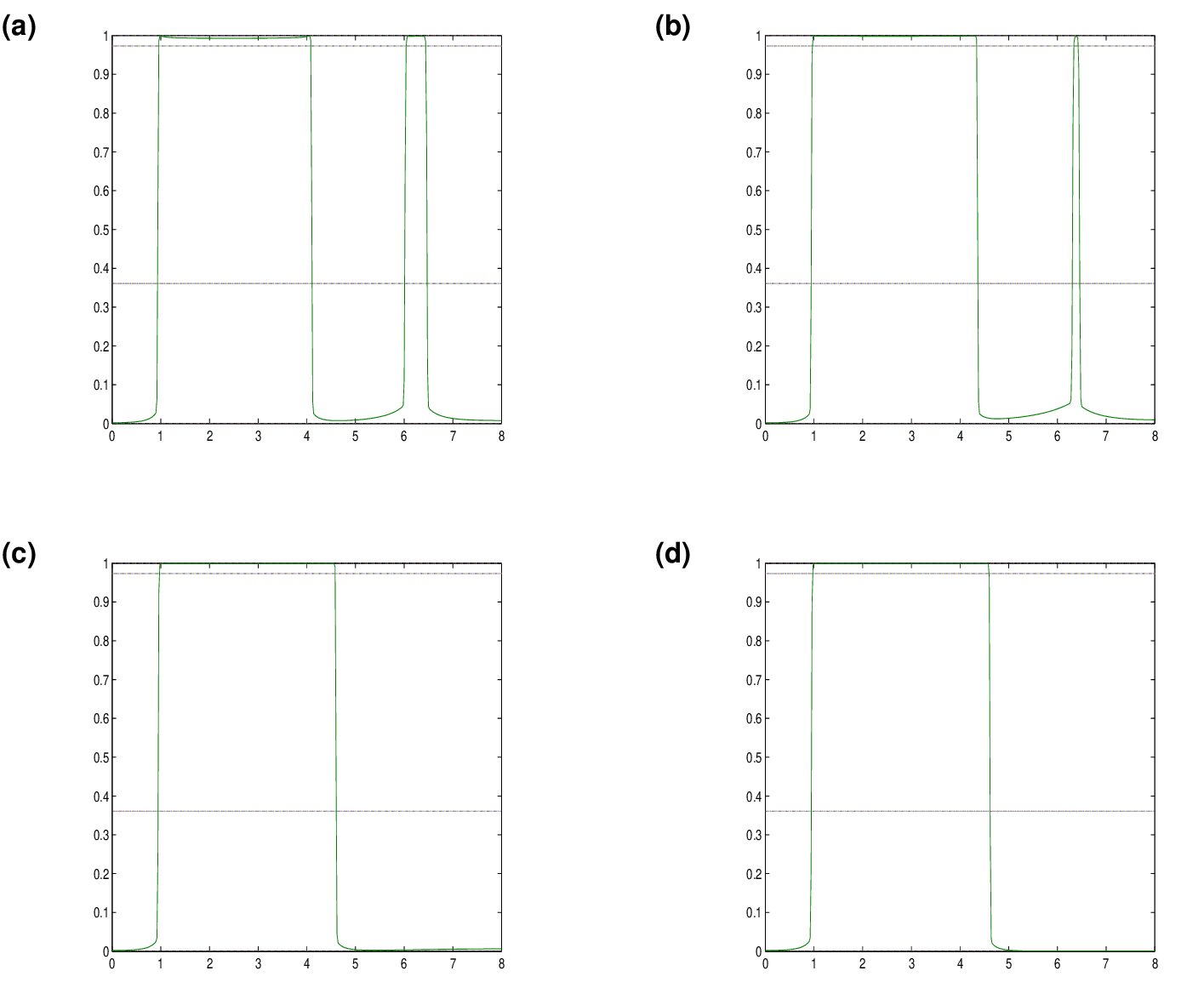}}

\caption{Evolution of an initial density profile having one large and one small plateau. Data shown at (a) $t=0$, (b) $t=8$, (c) $t=13$, (d) $t=18$. The parameter values are $\alpha=0.95$, $\chi_0=16$, $L=8$.}

\end{figure}

\subsection{Smooth, non-uniform steady states}

As mentioned below the statement of Theorem 4, our analytical results allow for the possibility that a solution of (\ref{rho})-(\ref{S})
with $\alpha>\frac{3}{4}$ might approach a smooth non-uniform steady state, avoiding $I_{\alpha}$, as $t\rightarrow\infty$. Some numerical evidence for this is presented in Figure 6, which depicts overlayed snapshots of a low-mass, high-adhesion solution converging to a bell-shaped steady state.

\begin{figure}

\centering

\resizebox{5in}{5in}{\includegraphics{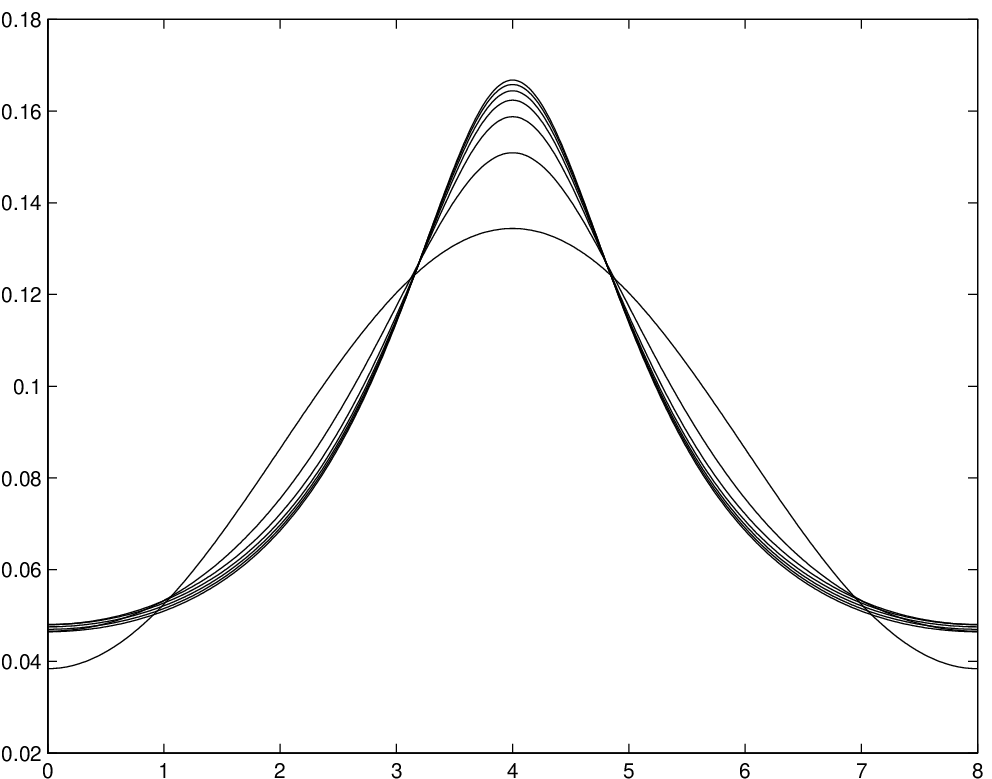}}

\caption{Evolution of a small initial density profile towards a non-uniform steady state below $I_{\alpha}=[0.361,0.973]$. Data shown at $t=0,2,4,6,8,10$ and $12$, such that the central maximum increases with time. The parameter values are $\alpha=0.95$, $\chi_0=16$, $L=8$.}

\end{figure}

\subsection{The question of critical mass}

For some chemotaxis models, such as the Keller-Segel equations in $\mathbb{R}^2$, to take a well-known example, there is a bifurcation phenomenon, such that solutions exist for all time (and also disperse) if the mass is below some critical value, while blow-up (formation of Dirac deltas) occurs in finite time otherwise \cite{pert}. We do not, however, expect our adhesion/chemotaxis model to exhibit quite this kind of bifurcation, since backward diffusion and volume filling have the effect of stabilising even very slender aggregations.

A numerical example of this is shown in Figure 7. Here, we took an initial datum with a thin high-density region, and a mass so small that the stability condition (\ref{boot}) is satisfied. Despite the uniform steady state being locally stable, the aggregation region appears to persist, such that the overlayed snapshots are visually indistinguishable. We would also expect to see analogous behaviour in an appropriate 2-d version of (\ref{rho})-(\ref{S}) in the high-adhesion regime.

A slightly different, but related, question one can ask is whether a low-mass aggregation satisfying (\ref{boot}) can be obtained by evolving an initial datum which lies below $I_{\alpha}$, but has, say, a narrow spike almost touching $\rho^{\flat}$. Despite numerous attempts, we have been unable to achieve this numerically; the spike always collapses almost immediately. Thus, there is evidence that (\ref{boot}) implies stability with respect to perturbations which remain below $\rho^{\flat}$.

\begin{figure}

\centering

\resizebox{5in}{5in}{\includegraphics{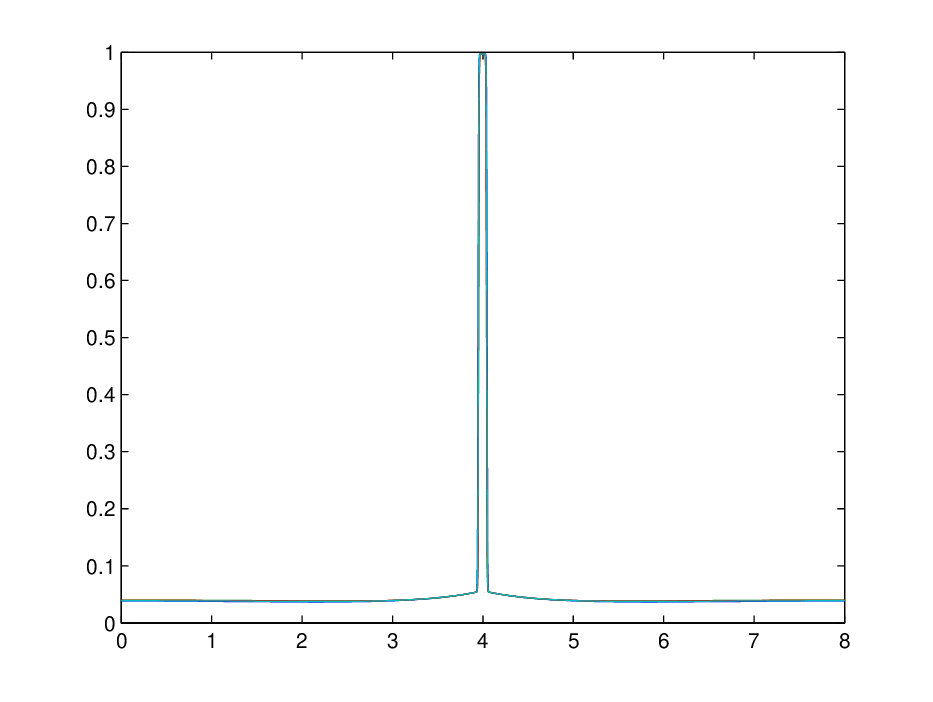}}

\caption{Overlayed snapshots of a low-mass aggregation at essentially steady state. Data shown at $t=0,0.2222,0.6667,1.1111$. The parameter values are $\alpha=0.95$, $\chi_0=16$, $L=8$.}

\end{figure}

\section{Stefan problems}

\subsection{Formulation}

In Section 4.1 we touched on an observation made in \cite{AS} for the special case $\chi_0=0$, namely that large-time plateau values in solutions of (\ref{walk}) seem to be essentially unique, for a given value of $\alpha$, and we noted that such uniqueness is also expected to hold for $\chi_0>0$, as a consequence of elliptic regularity. The observation of \cite{AS} subsequently led to the idea that a Stefan-problem framework, in which solutions are allowed to jump between unique plateau values $\rho_1(\alpha)$ and $\rho_2(\alpha)$ (as calculated in \cite{AS}), but are elsewhere smooth, might be an appropriate way of treating (\ref{rho}) as the limit of (\ref{walk}) in the high-adhesion regime, and it did indeed prove possible to develop an (at least partial) existence-and-uniqueness theory for such problems \cite{A}. Continuing in this vein, we will now proceed to write down a Stefan-problem formulation for (\ref{rho})-(\ref{S}) in the simplest possible case.

Imagine, then, that we are given a small initial density profile $\rho_0$ below the unstable region ({\em i.e.}, such that $\rho_0(x)<\rho^{\flat}~\forall x$), and imagine that we evolve this data {\em via} (\ref{rho})-(\ref{S}) until $\rho$ hits $\rho^{\flat}$ at some point $x_c$ and time $t_c$. The idea now is to continue the solution past $t_c$ by means of a singular, three-phase Stefan problem, whereby we introduce a high-density middle phase at $x_c$, which is initially and instantaneously a zero-width spike jumping up from $\rho_1$ to $\rho_2$ (and back down again), and which will subsequently fatten up as mass is drawn in from the low-density left-
and right-hand phases - thus, the width of the middle phase is strictly positive for $t>t_c$, and tends to zero as $t\searrow t_c$. In each of the left and right phases we impose the Dirichlet condition $\rho=\rho_1$ at the boundary with the middle phase, and in the middle phase we demand that $\rho=\rho_2$ at the left and right boundaries, which will be denoted by $s_l(t)$ and $s_r(t)$. The density in each phase evolves according to (\ref{rho}), the moving boundaries $s_l(t)$ and $s_r(t)$ evolve according to appropriate Rankine-Hugoniot conditions, and finally, since all of this is rather difficult to explain in words, we refer the reader to the simulation of Figure 8 (discussed below) for further clarification.

To be mathematically explicit, for $t>t_c$ we wish to solve in each phase the adhesion/chemotaxis equation
\be 
\frac{\p\rho}{\p t} = \frac{\p}{\p x}\left(D(\rho)\frac{\p\rho}{\p x}-\chi(\rho)\rho\frac{\p S}{\p x}\right),\label{rho_rep}
\ee
subject to the boundary conditions
\be
\frac{\partial\rho}{\partial x}=0\quad\textrm{at}\quad x=0\quad\textrm{and}\quad\rho=\rho_1\quad\textrm{at}\quad x=s_l^-(t)\label{bcl}
\ee
for the left-hand phase,
\be
\rho=\rho_2\quad\textrm{at}\quad x=s_l^+(t)\quad\textrm{and}\quad x=s_r^-(t)\label{bcm} 
\ee
for the middle (high-density) phase, and
\be
\rho=\rho_1\quad\textrm{at}\quad x=s_r^+(t)\quad\textrm{and}\quad\frac{\partial\rho}{\partial x}=0\quad\textrm{at}\quad x=L\label{bcr}
\ee
for the right-hand phase.

The chemoattractant concentration $S$ is obtained by globally solving the Neumann problem for 
\be
\Delta S = S-\rho\label{S_rep},
\ee
which is well-posed despite $\rho$ having two jump discontinuities, by elliptic regularity, as noted earlier.

The moving boundaries $s_l(t)$ and $s_r(t)$ are governed by the pair of Rankine-Hugoniot conditions
\be
\frac{ds_l}{dt} = -((D(\rho_1)\rho_x-\chi(\rho_1)\rho_1S_x)(s_l^-)-(D(\rho_2)\rho_x-\chi(\rho_2)\rho_2S_x)(s_l^+))/(\rho_1-\rho_2),\label{RHl}
\ee

\be
\frac{ds_r}{dt} = -((D(\rho_2)\rho_x-\chi(\rho_2)\rho_2S_x)(s_r^-)-(D(\rho_1)\rho_x-\chi(\rho_1)\rho_1S_x)(s_r^+))/(\rho_2-\rho_1),\label{RHr}
\ee
which guarantee local conservation of mass.

Unfortunately, equations (\ref{rho_rep})-(\ref{RHr}) have proved to be analytically intractable when subject to the singular initial condition $s_l(t_c)=s_r(t_c)$. In particular, we have been unable to prove the (plausible) conjecture that the Dirichlet condition $\rho=\rho_1$ at $s_l^-$ and $s_r^+$ {\em a priori} holds the density below $\rho^{\flat}$ in each of the low-density phases for some short time.

Nevertheless, we can at least attempt to solve these equations numerically, as we do below, whereby our attention will be focused on three questions:
\begin{enumerate}
\item Do sensible-looking solutions exist?
\item As a matter of principle, are solutions close (in some weak sense) to the oscillatory solutions obtained by discretising (\ref{rho})-(\ref{S}) directly?
\item Might it be more computationally efficient to solve the Stefan problem than to discretise (\ref{rho})-(\ref{S}) directly? 
\end{enumerate}

\subsection{The rescaled model; numerical solutions}

Following the approach of \cite{A}, we solve the $\rho$-equation in each given phase by rescaling the spatial coordinate so as to fix the relevant moving boundary (or boundaries). 

In the left-hand phase this results in
\be
\frac{\p\rho}{\p t} = \frac{1}{s_l^2}\frac{\p}{\p x}\left(D(\rho)\frac{\p\rho}{\p x}\right) + x\frac{\dot{s}_l}{s_l}\frac{\p\rho}{\p x}   -\frac{1}{s_l^2}\frac{\p}{\p x}\left(\chi(\rho)\rho\frac{\p S}{\p x}\right),\label{rho_left}
\ee
for $x\in[0,1]$, in the middle phase we get 
\be
\frac{\p\rho}{\p t} = \frac{1}{(s_r-s_l)^2}\frac{\p}{\p x}\left(D(\rho)\frac{\p\rho}{\p x}\right) + \frac{(\dot{s}_rx+(1-x)\dot{s}_l)}{(s_r-s_l)}\frac{\p\rho}{\p x} -\frac{1}{(s_r-s_l)^2}\frac{\p}{\p x}\left(\chi(\rho)\rho\frac{\p S}{\p x}\right),\label{rho_mid}
\ee
for $x\in[0,1]$, and in the right-hand phase
\be
\frac{\p\rho}{\p t} = \frac{1}{(L-s_r)^2}\frac{\p}{\p x}\left(D(\rho)\frac{\p\rho}{\p x}\right) + \frac{((1-x)\dot{s}_r)}{(L-s_r)}\frac{\p\rho}{\p x} -\frac{1}{(L-s_r)^2}\frac{\p}{\p x}\left(\chi(\rho)\rho\frac{\p S}{\p x}\right),\label{rho_right}
\ee
for $x\in[0,1]$.

The rescaled Rankine-Hugoniot conditions take the form
\be
\frac{ds_l}{dt} = -\left(\frac{1}{s_l}\left(D(\rho_1)\rho_x-\chi(\rho_1)\rho_1S_x\right)(1^-)-\frac{1}{(s_r-s_l)}\left(D(\rho_2)\rho_x-\chi(\rho_2)\rho_2S_x\right)(0^+)\right)/(\rho_1-\rho_2),\label{rescRHl}
\ee

\be
\frac{ds_r}{dt} = -\left(\frac{1}{(s_r-s_l)}\left(D(\rho_2)\rho_x-\chi(\rho_2)\rho_2S_x)(1^-\right)-\frac{1}{(L-s_r)}\left(D(\rho_1)\rho_x-\chi(\rho_1)\rho_1S_x\right)(0^+)\right)/(\rho_2-\rho_1).\label{rescRHr}
\ee

Each of equations (\ref{rho_left})-(\ref{rho_right}) is solved on a uniform grid in essentially the same way as in Section 4, while (\ref{rescRHl})-(\ref{rescRHr}) are solved (explicitly) using one-sided, second-order-accurate finite differences for the gradients ({\em e.g.}, $\rho_x(x_i)\approx\frac{1}{2h}(3\rho_i+\rho_{i-2}-4\rho_{i-1}))$. Since the three phases are generally of different physical lengths, this entails that the numerical approximation of $\rho$ lives on a globally-non-uniform grid. In order that we can nevertheless solve (\ref{S_rep}) conveniently, using the discrete Laplacian, $\rho$ is linearly interpolated onto a globally-uniform grid at each time step.

To obtain the simulation shown in Figure 8, we used the numerical method of Section 4 to evolve the initial data of Figure 4 (with $n=1200$ spatial points) until the solution hit $\rho=\rho^{\flat}$ at $t_c=0.8325$, and then continued the solution {\em via} the Stefan-problem algorithm just described, such that there are 100 spatial points in each of the three phases.

Note that the solution is nice and smooth away from the moving boundaries, and that the Dirichlet conditions at $s_l$ and $s_r$ hold the solution below $\rho^{\flat}$ for all time in the low-density phases. Also, since the gradient at $s_l^-$ and $s_r^+$ is large just after $t_c$, the middle phase gains mass very quickly for a short time; subsequent to this, there is a slow approach to the kind of weak steady-state solution of (\ref{rho})-(\ref{S}) seen in Figures 2-4.
 
In Figure 9 we overlay the simulations of Figures 4 and 8 in order to compare the Stefan-problem approach with that of direct discretisation. We see that, away from the central oscillatory region, there is always good agreement between the solutions, but that towards the middle of the domain there is a significant discrepancy shortly after $t_c$, due to the fact that, with the direct approach, the density has to push a considerable distance into $I_{\alpha}$ before oscillations set in (even when $n=1200$), thus creating a short delay. Also, the solution of Figure 4 unfortunately gains a little mass during the course of the simulation. Despite this, Figures 2-4 and 8, taken together, do seem to indicate that the Stefan problem is the correct (weak) limit of (\ref{walk}).

\begin{figure}

\centering

\resizebox{5in}{6in}{\includegraphics{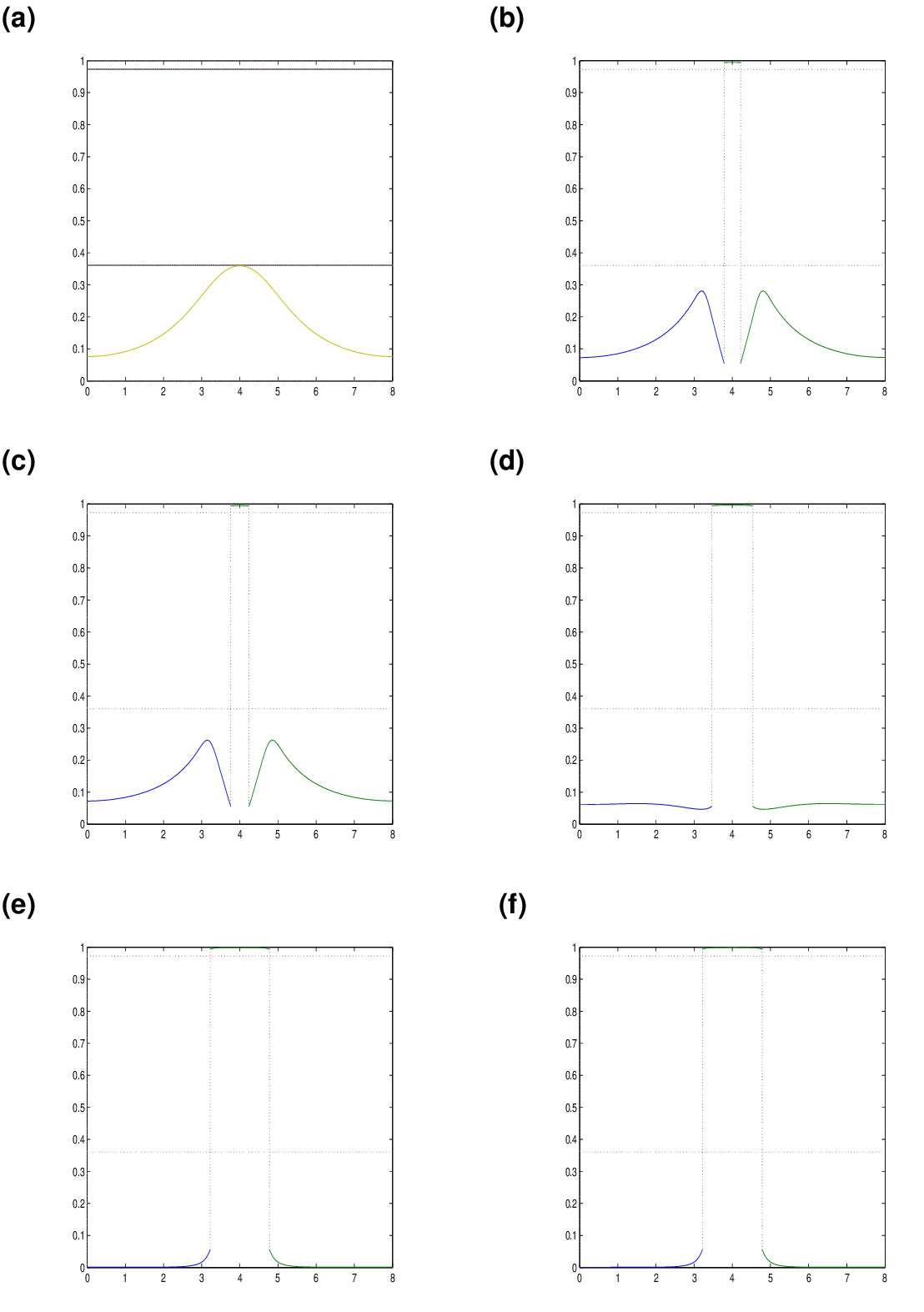}}

\caption{Here we show what happens when the $n=1200$ solution of Figure 4 is continued {\em via} the Stefan problem after hitting $I_{\alpha}$ at $t=0.8325$ (see (a)). In each phase there is a uniform mesh of 100 points. Snapshots (b)-(f) are taken at the same times as in Figure 4, namely (b) $t=1.0911$, (c) $t=1.1356$, (d) $t=2$, (e) $t=8.6778$, (f) $t=9.344$.}

\end{figure}

\begin{figure}

\centering

\resizebox{5in}{6in}{\includegraphics{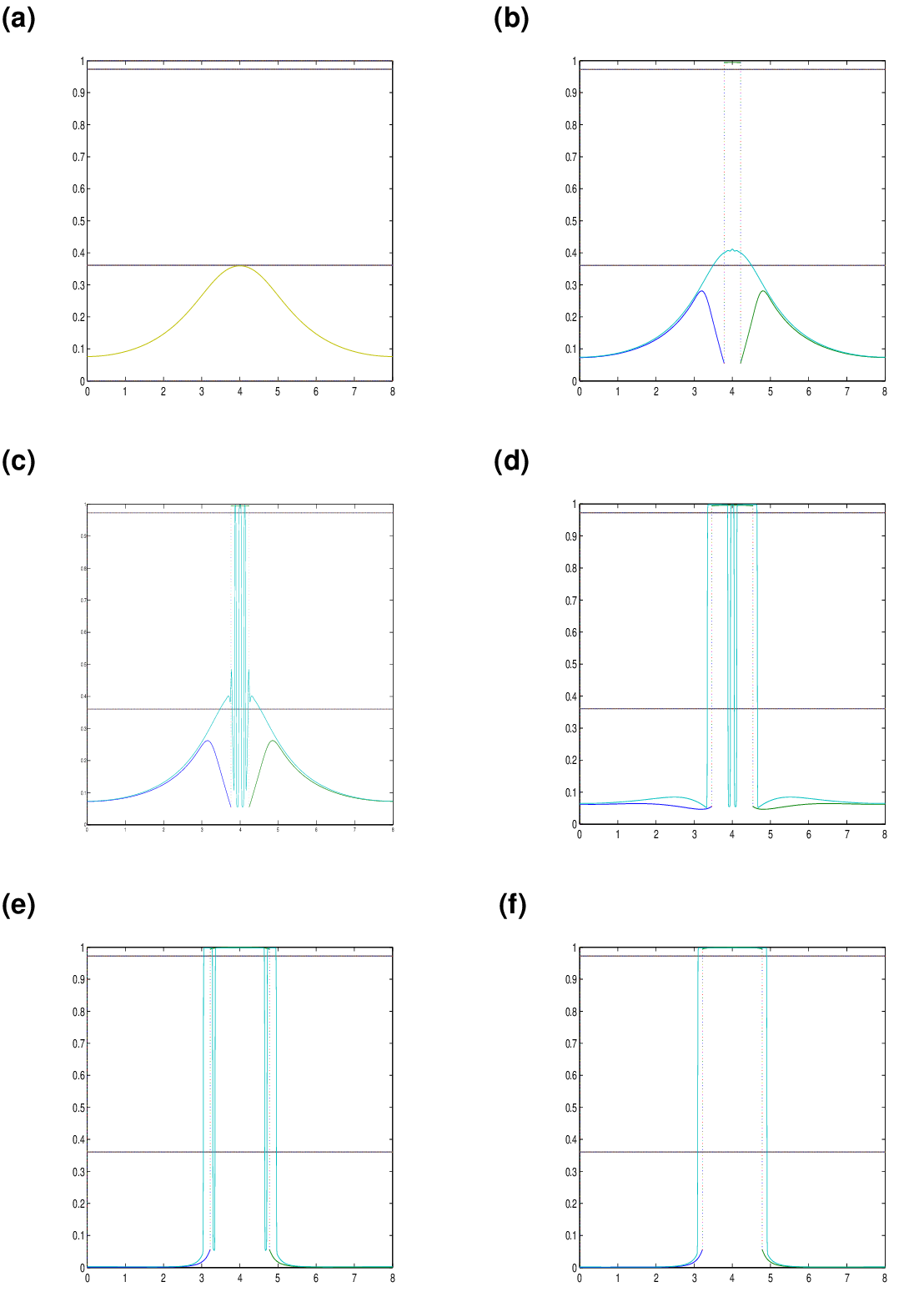}}

\caption{Initial data for the Stefan problem of Figure 7 (a), together with synchronised, overlayed snapshots from Figures 4 and 7 (b)-(f).}

\end{figure}

\section{Concluding remarks}
We have used a discrete random-walk model for cell adhesion and chemotaxis to generate sharp-edged cell aggregations from low-density initial data, and we have shown that a singular Stefan-problem description may be a fruitful way of approaching the ill-posed continuum-limit equations obtained in the high-adhesion regime.

One advantage of the Stefan-problem framework is that the microscopic oscillations seen in the underlying discrete model are avoided, and one can get convincing numerical solutions using a relatively coarse spatial grid. However, it should be noted that one {\em disadvantage} of the rather obvious numerical method we used for the three-phase problem (and which we certainly don't claim to be the best) is that, due to the spatial rescalings and the large initial gradients near the discontinuities, the parabolic and hyperbolic CFL conditions demand a very short time step until the high-density phase has attained a considerable thickness. 

Finally, although not simulated in this paper, one can of course imagine solutions of the three-phase Stefan problem in which chemotaxis is so strong that the solution in one of the low-density phases rises up to hit the unstable interval $I_{\alpha}$ once again. In that case, another spike should be inserted at the point of contact, and the solution continued {\em via} the appropriate five-phase Stefan problem, and so on, {\em ad infinitum}.

\section*{Acknowledgement}

The author wishes to thank the Wolfgang Pauli Institute, and in particular Christian Schmeiser, for supporting the completion of this work during a stay in Vienna.

\end{document}